\documentclass[a4paper,11pt]{amsart}

\usepackage{amssymb}
\usepackage{amsmath}
\usepackage{amsthm}
\usepackage[all]{xy}

\newcommand{\ii}{{\boldsymbol{i}}}
\newcommand{\half}{\tfrac12}

\newcommand{\CC}{\mathbb{C}}
\newcommand{\RR}{\mathbb{R}}
\newcommand{\ZZ}{\mathbb{Z}}

\newcommand{\X}{\mathsf{X}}
\newcommand{\Xv}{\X^\vee}
\newcommand{\Xvt}{\widetilde\X^\vee}
\newcommand{\Q}{\mathsf{Q}}
\newcommand{\Qv}{\Q^\vee}

\newcommand{\h}{\mathfrak{h}}
\newcommand{\s}{\mathfrak{s}}
\newcommand{\ttt}{\mathfrak{t}}

\newcommand{\e}{\varepsilon}
\newcommand{\ev}{\e^\vee}
\newcommand{\av}{\alpha^\vee}
\newcommand{\ov}{\omega^\vee}

\newcommand{\diag}{\operatorname{diag}}
\newcommand{\Hom}{\operatorname{Hom}}
\newcommand{\Img}{\operatorname{Im}}
\newcommand{\Ker}{\operatorname{Ker}}
\newcommand{\Gal}{\operatorname{Gal}}

\newcommand{\GL}{\mathrm{GL}}
\newcommand{\SL}{\mathrm{SL}}
\newcommand{\SO}{\mathrm{SO}}
\newcommand{\PSO}{\mathrm{PSO}}

\newcommand{\G}{\mathsf{G}}
\newcommand{\F}{\mathsf{F}}
\newcommand{\E}{\mathsf{E}}
\newcommand{\EV}{\mathsf{EV}}
\newcommand{\EVI}{\mathsf{EVI}}
\newcommand{\EVII}{\mathsf{EVII}}

\newcommand{\Gtil}{\widetilde{G}}
\newcommand{\Ztil}{\widetilde{Z}}

\newcommand{\uni}{\text{\upshape uni}}
\newcommand{\red}{\text{\upshape red}}
\newcommand{\sms}{{\rm ss}}
\newcommand{\ssc}{{\rm sc}}
\newcommand{\spl}{\text{\upshape s}}
\newcommand{\cmp}{\text{\upshape c}}

\newcommand{\sigmas}{{\sigma_\spl}}
\newcommand{\sigmac}{{\sigma_\cmp}}

\newcommand{\Ho}{\operatorname{H}}
\newcommand{\Zl}{\operatorname{Z}}
\newcommand{\Bd}{\operatorname{B}}

\newcommand{\Exp}{\mathcal{E}}
\newcommand{\Exptil}{\widetilde{\mathcal{E}}}

\newtheorem{theorem}{Theorem}
\newtheorem{lemma}[theorem]{Lemma}
\newtheorem{corollary}[theorem]{Corollary}

\theoremstyle{remark}
\newtheorem*{remark}{Remark}

\title[On the component group of a real reductive group]
{On the component group \\ of a real reductive group}

\author{Dmitry A. Timashev}

\address{Lomonosov Moscow State University, Faculty of Mechanics and
Mathematics, Department of Higher Algebra, 119991 Moscow, Russia}
\email{timashev@mccme.ru}

\thanks{This work was supported by the Russian Foundation
for Basic Research (grant 20-01-00091) and by the Ministry of Education and Science of the Russian Federation
in the framework of the program of the Moscow Center for Fundamental and Applied Mathematics (agreement 075-15-2019-1621).}



\begin{document}

\date{\today}

\begin{abstract}
For a connected linear algebraic group $G$ defined over~$\RR$, we compute the component group $\pi_0G(\RR)$ of the real Lie group $G(\RR)$ in terms of a maximal split torus $T_\spl\subseteq G$. In particular, we recover a theorem of Matsumoto (1964) that each connected component of $G(\RR)$ intersects~$T_\spl(\RR)$. We provide explicit elements of $T_\spl(\RR)$ which represent all connected components of~$G(\RR)$. The computation is based on structure results for real loci of algebraic groups and on methods of Galois cohomology.
\end{abstract}

\maketitle

\section*{Introduction}

Let $G$ be a linear algebraic group defined over the field of real numbers~$\RR$. The group of complex points $G(\CC)$ is a complex Lie group and the group of real points $G(\RR)$ is a real Lie group. Suppose that $G$ is connected (in the Zariski topology). It is well known that $G(\CC)$ is then connected (in the classical Hausdorff topology), but $G(\RR)$ is not necessarily connected: consider $G=\GL_n$ as an example. A natural problem is to compute the groups of connected components $\pi_0G(\RR)=G(\RR)/G(\RR)^\circ$, where $G(\RR)^\circ$ is the identity component of~$G(\RR)$.

In our previous paper joint with Mikhail Borovoi \cite{Galois-red} we computed $\pi_0G(\RR)$ in terms of real Galois cohomology. More precisely, there is an exact sequence
\begin{equation*}
1 \longrightarrow \pi_0G(\RR) \longrightarrow \Ho^1(\RR,\pi_1G(\CC)) \longrightarrow \Ho^1(\RR,\Gtil),
\end{equation*}
where $\pi_1G(\CC)$ is the fundamental group and $\Gtil$ is the universal cover of~$G(\CC)$, and $\Ho^1(\RR,{-})$ denotes the first (non-Abelian) group cohomology of $\Gal(\CC/\RR)$ acting naturally on a given group. Thus $\pi_0G(\RR)$ is the kernel of the natural map $\Ho^1(\RR,\pi_1G(\CC))\to\Ho^1(\RR,\Gtil)$.

In the same paper we developed combinatorial tools to compute the Galois cohomology set $\Ho^1(\RR,\Gtil)$. In fact, the computation reduces to the case of a simply connected semisimple group, treated in \cite{Galois-ss}. As for $\pi_1G(\CC)$, it is a finitely generated Abelian group and its Galois cohomology is well understood. In \cite{Galois-red} we described the map $\Ho^1(\RR,\pi_1G(\CC))\to\Ho^1(\RR,\Gtil)$ and derived a combinatorial description of $\pi_0G(\RR)$. In particular, it follows that $\pi_0G(\RR)$ is an elementary Abelian 2-group.

However, this description of the component group is somewhat abstract. For instance, it is not straightforward to produce from it explicit representatives of all connected components of~$G(\RR)$. On the other side, there is a classical result of Matsumoto \cite[Thm.\,1.2, Cor.]{Matsumoto} that each connected component of $G(\RR)$ intersects a given maximal $\RR$-split torus $T_\spl$ of $G$; see also \cite[Thm.\,14.4]{grp-red}. Hence $\pi_0G(\RR)$ is a quotient group of $\pi_0T_\spl(\RR)$, which explains geometrically, why it is an elementary Abelian 2-group.

Computation of real Galois cohomology of semisimple (or reductive) algebraic groups is close to classifying real forms or involutive automorphisms of these groups. There are two approaches to solving this kind of problems: the first one, going back to Kac, is via maximal anisotropic (or compact) tori and the second approach, going back to Satake, is via maximal split tori; see e.g.\ \cite[4.1, 4.4]{Lie-3}, \cite[5.1, 5.4]{Seminar-Lie}, \cite[26.5]{hom-emb}. In \cite{Galois-ss} and \cite{Galois-red} we used the ``Kac approach'' to compute $\Ho^1(\RR,G)$ based on a theorem of Borovoi \cite{Borovoi} reducing the computation to a maximal compact torus in~$G$. Here we use the ``Satake approach'' instead to obtain a more explicit and transparent description of $\pi_0G(\RR)$ in terms of a maximal split torus $T_\spl\subseteq G$.

The paper is organized as follows. After introducing notation and conventions and recalling basic definitions and facts in Section~\ref{s:prelim}, we briefly review real Galois cohomology in Section~\ref{s:Galois}. In Section~\ref{s:pi0} we consider the component group $\pi_0G(\RR)$ of an algebraic group $G$ over $\RR$ and obtain our main result (Theorem~\ref{t:pi0}). Examples are considered in Section~\ref{s:examples}.

\begin{remark}
Inspired by the first version of this paper on arXiv, Borovoi and Gabber gave a different proof of Theorem~\ref{t:pi0} in \cite{pi0-red}. Being shorter than our proof, it is based on Matsumoto's theorem and does not use Galois cohomology. On the contrary, our proof does not rely on Matsumoto's theorem and reestablishes a stronger version of this theorem instead.
\end{remark}

\section{Preliminaries}
\label{s:prelim}

\subsection{}

We fix a square root of $-1$ in $\CC$ and denote it by~$\ii$.

For an algebraic or Lie group $H$ we denote by $H^\circ$ its identity component and by $\h$ (the same lowercase German letter) its Lie algebra.

Throughout the paper, we denote by $G$ a connected linear algebraic group defined over $\RR$. We identify $G$ with its group of complex points $G(\CC)$ equipped with the antiregular involutive group automorphism $\sigma$ (\emph{real structure}) defining the action of $\Gal(\CC/\RR)$ on $G(\CC)$. (Here \emph{antiregular} means that $\sigma$ sends regular functions on Zariski open subsets of $G$ to complex conjugates of regular functions.) If $G$ is embedded in $\GL_n$ as a closed subgroup defined by polynomial equations with real coefficients in the matrix entries and inherits the standard real structure from $\GL_n$, then $\sigma$ is given by complex conjugation of the matrix entries. Note that $G(\RR)=G^\sigma$, the fixed point subgroup for~$\sigma$.

By a theorem of Mostow (see e.g.\ \cite[Thm.\,VIII.4.3]{Hochschild}) there is a Levi decomposition defined over $\RR$: $G=G_\uni\rtimes G_\red$, where $G_\uni$ is the unipotent radical of $G$ and $G_\red$ is reductive. The Levi subgroup $G_\red\subseteq G$ is defined up to conjugacy.

\subsection{}

Let $T$ be an algebraic torus of dimension~$n$. Then $T(\CC)\simeq(\CC^\times)^n$. Suppose that $T$ is equipped with a real structure~$\sigma$. Recall that $T$ is called \emph{$\RR$-split} (or just split) if $\sigma(t)=(\bar{t}_1,\dots,\bar{t}_n)$ and \emph{anisotropic} (or compact) if $\sigma(t)=(\bar{t}_1^{-1},\dots,\bar{t}_n^{-1})$ for any $t=(t_1,\dots,t_n)\in T(\CC)$. (Here a bar denotes complex conjugation.) If $T$ is split, then $T(\RR)\simeq(\RR^\times)^n$. If $T$ is anisotropic, then $T(\RR)\simeq(U_1)^n$ is a compact real torus, where $U_1$ denotes the unit circle in the complex plane. Any algebraic torus is equipped with a unique split real structure $\sigmas$ and with a unique anisotropic real structure~$\sigmac$.

A reductive connected algebraic group $G$ defined over $\RR$ is called \emph{anisotropic} (or compact) if $G$ does not contain $\RR$-split subtori. In this case $G(\RR)$ is a maximal compact real Lie subgroup of $G(\CC)$. It is connected by the polar decomposition theorem \cite[Thm.\,5.2.2]{Seminar-Lie}. An anisotropic real structure on a reductive group always exists and is unique up to conjugation; see e.g.\ \cite[5.2.3]{Seminar-Lie}. For any real structure $\sigma$ on $G$ there exists an anisotropic real structure $\sigmac$ \emph{compatible} with $\sigma$ in the sense that $\sigma$ and $\sigmac$ commute; see \cite[5.1.4]{Seminar-Lie} for the semisimple case and \cite[\S3]{Adams-Taibi} for the general reductive case. Then $\theta=\sigma\circ\sigmac=\sigmac\circ\sigma$ is a regular involutive group automorphism of $G$ commuting with both $\sigma$ and~$\sigmac$.

If $G=T$ is an algebraic torus, then we denote by $-\theta$ the involution given by $t\mapsto\theta(t)^{-1}$. There is an almost direct product decomposition $T=T_\cmp\cdot T_\spl$, where $T_\cmp=(T^\theta)^\circ$ is the unique maximal anisotropic subtorus and $T_\spl=(T^{-\theta})^\circ$ is the unique maximal split subtorus of~$T$. (Here \emph{almost direct product} means that the intersection $T_\cmp\cap T_\spl$ is finite, not necessarily trivial.) On the Lie algebra level, there is a decomposition $\ttt=\ttt_\cmp\oplus\ttt_\spl$ into the direct sum of $\theta$-eigenspaces with eigenvalues $1$ and~$-1$.

In the general case, there is a Cartan type decomposition
\begin{equation}\label{e:Cartan}
G^\sigmac=K\cdot T_\spl^\sigmac\cdot K,
\end{equation}
where $K=(G^\sigma)^\sigmac=(G^\theta)^\sigmac$ is a maximal compact subgroup both in $G^\sigma$ and in~$G^\theta$ and $T_\spl$ is a maximal split torus in~$G$; see e.g.\ \cite[Thm.\,V.6.7]{symm}. (In fact, the decomposition is established in [loc.~cit.] for semisimple $G$, but the reductive case is easily derived from the semisimple case.)

\subsection{}\label{ss:lattices}

Let $\X=\X(T)=\Hom(T,\CC^\times)$ denote the character lattice of~$T$. By identifying a character with its differential at the unity, we view $\X$ as a lattice in~$\ttt^*$. The dual lattice $\Xv=\Xv(T)=\Hom(\CC^\times,T)$ is the cocharacter lattice; it can be regarded as a lattice in~$\ttt$.

We put $\Xv_\cmp=\Xv\cap\ttt_\cmp$ and $\Xv_\spl=\Xv\cap\ttt_\spl$; these are the cocharacter lattices of $T_\cmp$ and~$T_\spl$, respectively. The  involutions $\sigma$ and $\theta$ act on $\Xv$ in a natural way, so that $\theta$ acts on $\Xv_\cmp$ identically and on $\Xv_\spl$ as $-1$, while $\sigma$ acts as $-\theta$.

We also denote by $\Xvt_\cmp$ and $\Xvt_\spl$ the images of $\Xv$ under the projection maps from $\ttt$ to $\ttt_\cmp$ and $\ttt_\spl$ given by
\begin{equation*}
\nu\mapsto\nu_\cmp=\half(\nu+\theta(\nu))\quad\text{and}\quad\nu\mapsto\nu_\spl=\half(\nu-\theta(\nu)),
\end{equation*}
respectively. Note that
\begin{equation*}
\Xv_\cmp\subseteq\Xvt_\cmp\subseteq\half\Xv_\cmp\quad\text{and}\quad\Xv_\spl\subseteq\Xvt_\spl\subseteq\half\Xv_\spl.
\end{equation*}
Similar notation $\Lambda_\cmp$, $\Lambda_\spl$, $\widetilde\Lambda_\cmp$, $\widetilde\Lambda_\spl$ will be used for any $\theta$-stable sublattice $\Lambda\subset\Xv$.

\section{Galois cohomology over real numbers}
\label{s:Galois}

\subsection{}

Our basic reference on Galois cohomology is Serre's book \cite{Galois}. We consider non-Abelian cohomology of the Galois group $\Gal(\CC/\RR)$ of order~2 with coefficients in linear algebraic groups over $\RR$ or even in arbitrary groups on which $\Gal(\CC/\RR)$ acts by group automorphisms. Such an action of $\Gal(\CC/\RR)$ on a group $A$ is defined by an involutive group automorphism $\sigma$ of~$A$.

A \emph{$1$-cocycle} with coefficients in $A$ is an element $z\in A$ such that $z\cdot\sigma(z)=1$. The set $\Zl^1(\RR,A)$ of all 1-cocycles is endowed with an $A$-action by twisted conjugation:
\begin{equation*}
z\overset{a}\longmapsto a\cdot z\cdot\sigma(a)^{-1}.
\end{equation*}
The orbits $[z]$ of this action are the cohomology classes and the orbit set
\begin{equation*}
\Ho^1(\RR,A)=\Zl^1(\RR,A)/A
\end{equation*}
is the \emph{first cohomology set} with coefficients in~$A$ or, for brevity, the first Galois cohomology set of~$A$. It is a pointed set, the base point being $[1]$, the class of the trivial 1-cocycle, also denoted as $\Bd^1(\RR,A)$ (the set of 1-cobound\-aries).

For a short exact sequence of groups
\begin{equation}\label{e:short}
1 \longrightarrow A \longrightarrow B \longrightarrow C \longrightarrow 1
\end{equation}
equipped with compatible actions of~$\Gal(\CC/\RR)$ there is a long exact cohomology sequence
\begin{equation}\label{e:long-cohom}
1 \longrightarrow A^\sigma \longrightarrow B^\sigma \longrightarrow C^\sigma \longrightarrow
\Ho^1(\RR,A) \longrightarrow \Ho^1(\RR,B) \longrightarrow \Ho^1(\RR,C);
\end{equation}
see \cite[I.5.5]{Galois}. (Here the same letter $\sigma$ denotes an involution of $B$ and the induced involutions of its subgroup $A$ and the quotient group~$C$.) The maps in \eqref{e:long-cohom} are obvious, except for the fourth one, given by
\begin{equation}\label{e:connect}
c\longmapsto[b^{-1}\sigma(b)],\qquad\forall c\in C^\sigma,
\end{equation}
where $b\in B$ is any element in the preimage of~$c$, so that $b^{-1}\sigma(b)\in A$.
Exactness is understood in the sense of pointed sets. In particular, triviality of the image of a map does not readily imply injectivity of the subsequent map. However the fibers of the maps in \eqref{e:long-cohom} can be identified with the kernels of the maps in similar long exact cohomology sequences, corresponding to the short exact sequences \eqref{e:short} with twisted actions of $\Gal(\CC/\RR)$; see \cite[I.5.5]{Galois} for details.

\subsection{}

Now let $G$ be a connected linear algebraic group defined over~$\RR$. We have a long exact cohomology sequence
\begin{equation}\label{e:cohom-Levi}
\dots \longrightarrow \Ho^1(\RR,G_\uni) \longrightarrow \Ho^1(\RR,G) \longrightarrow \Ho^1(\RR,G/G_\uni).
\end{equation}
Due to Galois cohomology vanishing for unipotent groups (see \cite[Lemma 6.2(ii)]{Galois-red} for an elementary proof over~$\RR$) the last map in \eqref{e:cohom-Levi} is injective. It is also surjective, because a Levi subgroup $G_\red\subseteq G$ maps isomorphically onto $G/G_\uni$. Hence $\Ho^1(\RR,G)\simeq\Ho^1(\RR,G_\red)$ and the computation of Galois cohomology boils down to the reductive case. See Theorem~6.5(i) of \cite{Galois-red} and its proof for details.

For reductive $G$, the computation reduces to a maximal compact subgroup~$G^\sigmac$.

\begin{lemma}[{\cite[Lemma 3.2]{R-orb-symm}}]\label{l:unitary}
Suppose $G$ is a connected reductive group defined over $\RR$. Each Galois cohomology class in $\Ho^1(\RR,G)$ is represented by a 1-cocycle $z\in G^\sigmac$. Two cocycles $z,z'\in G^\sigmac$ are in the same cohomology class if and only if $z'=g\cdot z\cdot\sigma(g)^{-1}$ for some $g\in G^\sigmac$.
\end{lemma}

\section{The component group of a real algebraic group}
\label{s:pi0}

\subsection{}

In computing the component group $\pi_0G(\RR)$ we follow the strategy of \cite{Galois-red}. First we note that the Levi decomposition of $G$ yields a similar decomposition of the real locus: $G(\RR)=G_\uni(\RR)\rtimes G_\red(\RR)$. Since $G_\uni(\RR)$ is connected (this is proved using the exponential map; see e.g.\ \cite[Lemma 6.2(i)]{Galois-red}), we have $\pi_0G(\RR)=\pi_0G_\red(\RR)$. Thus the problem is reduced to the case of a reductive group.

Assume now that $G$ is reductive. There is an almost direct product decomposition $G=G_\sms\cdot S$, where $G_\sms$ is a connected semisimple group and $S$ is a torus. In fact, $G_\sms=[G,G]$ is the derived subgroup and $S=Z(G)^\circ$ is the connected center of~$G$. Let $G_\ssc$ denote the universal cover of~$G_\sms$; it is a simply connected semisimple group. Then $\Gtil=G_\ssc\times\s$ is the universal cover of~$G$ and we have a natural exact sequence
\begin{equation}\label{e:uni-cover}
1 \longrightarrow \pi_1{G} \overset{i}{\longrightarrow} \Gtil=G_\ssc\times\s \overset{j}{\longrightarrow} G=G_\sms\cdot S \longrightarrow 1,
\end{equation}
where $\pi_1G$ is the fundamental group of $G(\CC)$. The covering map $j$ is given by the formula
\begin{equation*}
j(\tilde{g})=j_\ssc(g_\ssc)\cdot\Exp(y)\quad\text{for}\quad\tilde{g}=(g_\ssc,y),\ g_\ssc\in G_\ssc,\ y\in\s,
\end{equation*}
where $j_\ssc:G_\ssc\to G_\sms$ is the covering map of the derived subgroup and $\Exp:\s\to S$ is given by the formula
\begin{equation}\label{e:Exp}
\Exp(y)=\exp{2\pi y}.
\end{equation}
Note that $\Gtil$ can be regarded as an algebraic group defined over $\RR$ (if we interpret $\s$ as an Abelian unipotent group) and the map $j$ commutes with the action of $\Gal(\CC/\RR)$, but $j$ is a homomorphism of complex Lie groups, \emph{not} of algebraic groups.

The fundamental group $\pi_1G$ also carries on a natural action of $\Gal(\CC/\RR)$ induced from the action on~$G$ and the map $i$ commutes with this action. The group $\Ztil=\Img(i)=\Ker(j)$ is a discrete central subgroup of~$\Gtil$.

Let $T$ be a maximal torus in $G$ defined over~$\RR$. There is an almost direct product decomposition $T=T_\sms\cdot S$, where $T_\sms=T\cap G_\sms$ is a maximal torus in~$G_\sms$. Then $T_\ssc=j_\ssc^{-1}(T_\sms)$ is a maximal torus in $G_\ssc$ and $\Xv(T_\ssc)=\Qv$ is the coroot lattice of~$G$.

There is a commutative diagram with exact rows and columns:
\begin{equation*}
\xymatrix{
         &              &  0  & 0  &   \\
1 \ar[r] & \Ztil \ar[r] &    T_\ssc\times\s \ar[r]^-j \ar[u]    & T \ar[r] \ar[u] & 1 \\
         &     0 \ar[r] & \ttt_\ssc\oplus\s \ar[r]^-\sim \ar[u]^\Exptil & \ttt \ar[r] \ar[u]^\Exp & 0 \\
         &     0 \ar[r] & \ii\Qv \ar[r] \ar[u] & \ii\Xv \ar[u] &   \\
         &              &  0 \ar[u] & 0 \ar[u] &   \\
}\end{equation*}
Here the map $\Exp:\ttt\to T$ is given by the formula \eqref{e:Exp} and $\Exptil(x,y)=(\Exp_\ssc(x),y)$ for $x\in\ttt_\ssc$, $y\in\s$, where $\Exp_\ssc:\ttt_\ssc\to T_\ssc$ is the scaled exponential map similar to~\eqref{e:Exp}. It follows that
\begin{equation*}
\pi_1G\simeq\Ztil=\Exptil(\ii\Xv)\simeq\ii\Xv\!/\,\ii\Qv
\end{equation*}
is a finitely generated Abelian group.

\subsection{}

The Galois cohomology exact sequence \eqref{e:long-cohom} of the short exact sequence \eqref{e:uni-cover} reads as
\begin{equation*}
\cdots \longrightarrow \Gtil(\RR) \longrightarrow G(\RR) \longrightarrow \Ho^1(\RR,\pi_1{G}) \longrightarrow \Ho^1(\RR,\Gtil) \longrightarrow \cdots
\end{equation*}
By a theorem of \'E.~Cartan (see e.g.\ \cite[Cor.\,4.7]{grp-red-add} or \cite[Thm.\,4.2.2]{Lie-3}) $G_\ssc(\RR)$ is a connected Lie group, whence $\Gtil(\RR)=G_\ssc(\RR)\times\s(\RR)$ is connected, too. By dimension reasons, $\Gtil(\RR)$ maps onto $G(\RR)^\circ$. Since the Galois cohomology of the Abelian unipotent group $\s$ is trivial, we have
\begin{equation}\label{e:H1}
\Ho^1(\RR,\Gtil) \simeq \Ho^1(\RR,G_\ssc) \times \Ho^1(\RR,\s) \simeq \Ho^1(\RR,G_\ssc).
\end{equation}
Since $\pi_1G\simeq\ii\Xv/\ii\Qv$, we finally get an exact sequence
\begin{equation}\label{e:exact}
1 \longrightarrow \pi_0G(\RR) \overset\delta\longrightarrow \Ho^1(\RR,\ii\Xv\!/\,\ii\Qv) \overset{i_*}\longrightarrow \Ho^1(\RR,G_\ssc).
\end{equation}
The map $\delta$ sends a connected component $gG(\RR)^\circ$ of $G(\RR)$ to the cohomology class $[\ii\nu+\ii\Qv]$ such that $\tilde{z}:=\Exptil(\ii\nu)=\tilde{g}^{-1}\sigma(\tilde{g})\in\Zl^1(\RR,\Ztil)$, where $\tilde{g}\in\Gtil$ is such that $j(\tilde{g})=g$.
The map $i_*$ sends a cohomology class $[\ii\nu+\ii\Qv]$ to the cohomology class of $z_\ssc=\Exp_\ssc(\ii\nu_\ssc)\in G_\ssc$, where $\nu_\ssc$ is the projection of $\nu\in\ttt=\ttt_\ssc\oplus\s$ to~$\ttt_\ssc$.

\subsection{}

Now we compute the first cohomology group of $\ii\Xv\!/\,\ii\Qv$. (Note that it is indeed a group, because $\ii\Xv\!/\,\ii\Qv$ is Abelian.)

\begin{lemma}\label{l:ZB1(pi1)}
For $\nu\in\Xv$ we have: $\ii\nu+\ii\Qv\in\Zl^1(\RR,\ii\Xv\!/\,\ii\Qv)$ if and only if $\nu_\cmp\in\half\Qv_\cmp$. We also have: $\ii\nu+\ii\Qv\in\Bd^1(\RR,\ii\Xv\!/\,\ii\Qv)$ if and only if $\nu\in2\Xvt_\spl+\Qv$.
\end{lemma}

\begin{proof}
The action of $\sigma$ on $\ii\Xv$ coincides with the action of~$\theta$. Hence
\begin{equation*}
\Zl^1(\RR,\ii\Xv\!/\,\ii\Qv)=(\ii\Xv\!/\,\ii\Qv)^{-\theta},
\end{equation*}
the fixed point set of~$-\theta$. Thus $\ii\nu+\ii\Qv$ is a 1-cocycle if and only if $\theta(\nu)=-\nu$ modulo $\Qv$, i.e., $\nu+\theta(\nu)=2\nu_\cmp\in\Qv$.

The set of 1-coboundaries $\Bd^1(\RR,\ii\Xv\!/\,\ii\Qv)$ is the image of $\ii\Xv\!/\,\ii\Qv$ under the induced action of $1-\theta$. Thus each coboundary coset is of the form $\ii\nu+\ii\Qv$ with $\nu=\mu-\theta(\mu)=2\mu_\spl$ for some $\mu\in\Xv$. The claim follows.
\end{proof}

\begin{corollary}\label{c:H1(pi1)}
$\Ho^1(\RR,\ii\Xv\!/\,\ii\Qv)\simeq\Xv\!\cap(\Xvt_\spl+\half\Qv_\cmp)/(2\Xvt_\spl+\Qv)$.
\end{corollary}

\subsection{}

Assume from now on that $T_\spl$ is a \emph{maximal} $\RR$-split torus in $G$ and $T$ is a maximal torus containing~$T_\spl$. We describe the kernel of the map $i_*$ in \eqref{e:exact}.

\begin{lemma}
A 1-cocycle $\tilde{z}=\Exptil(\ii\nu)\in\Zl^1(\RR,\Ztil)$ is a 1-coboundary in $\Gtil$ if and only if $z_\ssc=\Exp_\ssc(\ii\nu_\ssc)\in(T_\ssc)_\spl$.
\end{lemma}

\begin{proof}
By \eqref{e:H1} we have $\tilde{z}\in\Bd^1(\RR,\Gtil)$ if and only if $z_\ssc=\Exp_\ssc(\ii\nu_\ssc)\in\Bd^1(\RR,G_\ssc)$. Since $z_\ssc$ is a central element of $G_\ssc$, it has finite order, hence $z_\ssc\in G_\ssc^\sigmac$. By Lemma~\ref{l:unitary}, $z_\ssc=g_\ssc^{-1}\sigma(g_\ssc)$ for some $g_\ssc\in G_\ssc^\sigmac$.

By \eqref{e:Cartan} there is a decomposition $g_\ssc=k_1 t_\ssc k_2$, where $k_1,k_2\in K$ and $t_\ssc\in(T_\ssc)_\spl^\sigmac$. Then $z_\ssc=k_2^{-1}t_\ssc^{-2}k_2=t_\ssc^{-2}$, because $z_\ssc$ is central. Hence $z_\ssc\in(T_\ssc)_\spl$.

Conversely, for any $s\in(T_\ssc)_\spl^\sigmac$ there exists $t_\ssc\in(T_\ssc)_\spl^\sigmac$ such that $s=t_\ssc^{-2}=t_\ssc^{-1}\sigma(t_\ssc)$. Hence $s$ is a 1-coboundary in $(T_\ssc)_\spl$ and more so in~$G_\ssc$.
\end{proof}

\begin{corollary}\label{c:pi0}
$\Ker(i_*)\simeq(\Xv_\spl+\Qv)/(2\Xvt_\spl+\Qv)\simeq\Xv_\spl/(2\Xvt_\spl+\Qv_\spl)$.
\end{corollary}

\begin{proof}
With the above notation, $[\ii\nu+\Qv]\in\Ker(i_*)$ if and only if ${z_\ssc\in(T_\ssc)_\spl}$. Since $\nu$ is defined modulo $\Qv\subset\ttt_\ssc$ and $\nu_\cmp\in\ttt_\ssc$ by Lemma~\ref{l:ZB1(pi1)}, we may assume that $\nu_\ssc\in(\ttt_\ssc)_\spl$ and $\nu\in\ttt_\spl$, whence $\nu\in\Xv_\spl$. We conclude by Corollary~\ref{c:H1(pi1)}.
\end{proof}

\subsection{}

We are now ready to state and prove our main result.

\begin{theorem}\label{t:pi0}
Let $G$ be a connected linear algebraic group defined over~$\RR$. Choose a maximal $\RR$-split torus $T_\spl\subseteq G$, and a maximal torus $T\subseteq G$ defined over $\RR$ and containing~$T_\spl$. With the notation of Subsection~\ref{ss:lattices},
\begin{equation}\label{e:pi0}
\pi_0G(\RR)\simeq\Xv_\spl/(2\Xvt_\spl+\Qv_\spl),
\end{equation}
where $\Xv$ is the cocharacter lattice of $T$ and $\Qv$ is the coroot lattice of $G/G_\uni$ with respect to (the image of)~$T$.
The connected component of $G(\RR)$ corresponding to the coset of $\nu\in\Xv_\spl$ is represented by the element
\begin{equation*}
t=\Exp(\ii\nu/2)=\exp\pi\ii\nu
\end{equation*}
of order dividing 2 in~$T_\spl(\RR)$.
\end{theorem}

\begin{proof}
Choosing a Levi subgroup $G_\red\subseteq G$ containing $T$ reduces the proof to the case of reductive~$G$. In this case, $\pi_0G(\RR)$ is given by the formula \eqref{e:pi0} due to the exact sequence \eqref{e:exact} and Corollary~\ref{c:pi0}.

For any $\nu\in\Xv_\spl$, consider the respective 1-cocycle $\tilde{z}=\Exptil(\ii\nu)\in\Zl^1(\RR,\Ztil)$. We have $\tilde{z}=\tilde{t}^{-2}=\tilde{t}^{-1}\sigma(\tilde{t})$ for $\tilde{t}=\Exptil(-\ii\nu/2)$. Then $t=j(\tilde{t})$ and the map $\delta$ of \eqref{e:exact} sends the connected component of $G(\RR)$ containing $t$ to the cohomology class of $\ii\nu+\ii\Qv$ in $\Ho^1(\RR,\ii\Xv\!/\,\ii\Qv)$. This completes the proof.
\end{proof}

\section{Examples}
\label{s:examples}

\subsection{}

Suppose that $G=T$ is an algebraic torus defined over~$\RR$. In this case we have
\begin{equation*}
\pi_0T(\RR)\simeq\Xv_\spl/2\Xvt_\spl\simeq T_\spl^{(2)}/(T_\spl\cap T_\cmp),
\end{equation*}
where $T_\spl^{(2)}$ is the elementary Abelian 2-subgroup of elements of order $\le2$ in~$T_\spl$. The second isomorphism sends the coset of $\nu\in\Xv_\spl$ to the coset of $t=\exp\pi\ii\nu\in T_\spl^{(2)}$. Note that $\nu\in2\Xvt_\spl$ if and only if $\nu-\mu\in2\Xv$ for some $\mu\in\ttt_\cmp$, i.e., $t=\exp\pi\ii\mu\in T_\spl\cap T_\cmp$.

\subsection{}\label{GL}

Let $G=\GL_n$. We may choose a maximal torus $T=T_\spl$ consisting of diagonal matrices. The weights $\e_1,\dots,\e_n$ of the tautological $n$-dimensional representation of $\GL_n$ constitute a basis of $\X$. Let $\{\ev_1,\dots,\ev_n\}$ be the dual basis of~$\Xv$. The roots of $\GL_n$ with respect to $T$ are $\e_i-\e_j$ and the coroots are $\ev_i-\ev_j$ over all $i\ne j$. The coroot lattice is
\begin{equation*}
\Qv=\{k_1\ev_1+\dots+k_n\ev_n\mid k_i\in\ZZ,\ k_1+\dots+k_n=0\}.
\end{equation*}
Since $\Xv=\Xv_\spl=\Xvt_\spl$, we have
\begin{equation*}
2\Xvt_\spl+\Qv_\spl=\{k_1\ev_1+\dots+k_n\ev_n\mid k_i\in\ZZ,\ k_1+\dots+k_n\text{ is even}\}
\end{equation*}
and therefore $\pi_0\GL_n(\RR)$ is of order~2, as expected. The only nonzero element of $\Xv_\spl/(2\Xvt_\spl+\Qv_\spl)$ is the coset of $\ev_1$ and the respective element $t\in T_\spl$ representing the non-unity connected component $\GL_n(\RR)\setminus\GL_n(\RR)^\circ$ is
\begin{equation*}
t=\exp\pi\ii\ev_1=\diag(-1,1,\dots,1).
\end{equation*}

\subsection{}

Let $G=\SO_{p,q}$ with $0<p\le q$. Then $\dim T=\ell:=\lfloor n/2\rfloor$, where $n=p+q$, and $\dim T_\spl=p$. We may assume that $\SO_{p,q}$ preserves the indefinite inner product
\begin{equation*}
(x,y)=x_1y_n+\dots+x_py_{q+1}+x_{p+1}y_{p+1}+\dots+x_qy_q+x_{q+1}y_p+\dots+x_ny_1.
\end{equation*}
The standard real structure on $\SO_{p,q}$ is given by complex conjugation of the matrix entries. We may choose $T_\spl$ consisting of the matrices of the form
\begin{equation*}
t=\diag(t_1,\dots,t_p,1,\dots,1,t_p^{-1},\dots,t_1^{-1})
\end{equation*}
and $T_\cmp$ consisting of the matrices of the form
\begin{equation*}
t=\diag(\underbrace{1,\dots,1}_p,U_{p+1},\dots,U_{\ell},\underbrace{1,\dots,1}_{p\text{ or }p+1}),\quad\text{where}\quad
U_i=\begin{pmatrix}
 \frac{t_i+\text{\rlap{$t_i^{-1}$}}\hphantom{t_i}}2      & \frac{t_i^{-1}-t_i}{2\ii} \\
 \frac{t_i-\text{\rlap{$t_i^{-1}$}}\hphantom{t_i}}{2\ii} & \frac{t_i+\text{\rlap{$t_i^{-1}$}}\hphantom{t_i}}2
\end{pmatrix}.
\end{equation*}

Let $\pm\e_1,\dots,\pm\e_\ell$ be the nonzero weights of the tautological $n$-dimensional representation of $\SO_n$. With the above notation, $\e_i(t)=t_i$. Then $\{\e_1,\dots,\e_\ell\}$ is a basis of~$\X$. Let $\{\ev_1,\dots,\ev_\ell\}$ be the dual basis of~$\Xv$. The roots of $\SO_{p,q}$ with respect to $T$ are $\pm\e_i\pm\e_j$ over all $i\ne j$ and, for odd~$n$, $\pm\e_i$ ($i=1,\dots,\ell$). The coroots are $\pm\ev_i\!\!\pm\ev_j$ and, for odd~$n$, $\pm2\ev_i$. The coroot lattice is
\begin{equation*}
\Qv=\{k_1\ev_1+\dots+k_\ell\ev_\ell\mid k_i\in\ZZ,\ k_1+\dots+k_\ell\text{ is even}\}.
\end{equation*}
The action of $\theta$ on $\Xv$ is given by the formula
\begin{equation*}
\theta(\ev_i)=
\begin{cases}
            -\ev_i, & i=1,\dots,p,   \\
 \hphantom{-}\ev_i, & i=p+1,\dots,\ell.
\end{cases}
\end{equation*}

The lattice $\Xv_\spl=\Xvt_\spl$ is spanned by $\ev_1,\dots,\ev_p$ and
\begin{equation*}
\Qv_\spl=\{k_1\ev_1+\dots+k_p\ev_p\mid k_i\in\ZZ,\ k_1+\dots+k_p\text{ is even}\}\supseteq2\Xvt_\spl.
\end{equation*}
We conclude that $\pi_0\SO_{p,q}(\RR)$ is a group of order~2 (a well-known fact). The only nonzero element of $\Xv_\spl/(2\Xvt_\spl+\Qv_\spl)$ is the coset of $\ev_1$ and the respective element $t\in T_\spl$ representing the non-unity connected component $\SO_{p,q}(\RR)\setminus\SO_{p,q}(\RR)^\circ$ is
\begin{equation*}
t=\exp\pi\ii\ev_1=\diag(-1,1,\dots,1,-1).
\end{equation*}

\subsection{}

We modify the previous example by taking $G=\PSO_{p,q}=\SO_{p,q}/\{\pm E\}$ with even $n=p+q=2\ell$. We use the notation from the previous subsection. The cocharacter lattice here is twice bigger:
\begin{equation*}
\Xv=\{k_1\ev_1+\dots+k_\ell\ev_\ell\mid k_i\in\half\ZZ,\ k_i-k_j\in\ZZ,\ \forall i,j=1,\dots,\ell\}.
\end{equation*}

If $p<q$, then $\Xv_\spl$ is spanned by $\ev_1,\dots,\ev_p$ and
\begin{equation*}
\Xvt_\spl=\{k_1\ev_1+\dots+k_p\ev_p\mid k_i\in\half\ZZ,\ k_i-k_j\in\ZZ,\ \forall i,j=1,\dots,p\}.
\end{equation*}
For even $p$ we get $2\Xvt_\spl\subseteq\Qv_\spl$ and therefore $\pi_0\PSO_{p,q}(\RR)$ is of order~2, with the non-unity component represented by the same matrix $t$ as in the previous example. For odd $p$ we get $2\Xvt_\spl+\Qv_\spl=\Xv_\spl$ and therefore $\PSO_{p,q}(\RR)$ is connected in this case.

If $p=q=\ell$ (the split case), then $\Xv_\spl=\Xvt_\spl=\Xv$. If $\ell$ is odd, then
\begin{equation*}
2\Xvt_\spl+\Qv_\spl=\{k_1\ev_1+\dots+k_\ell\ev_\ell\mid k_i\in\ZZ\}
\end{equation*}
and $\pi_0\PSO_{\ell,\ell}(\RR)$ is of order~2. The only nonzero element of $\Xv_\spl/(2\Xvt_\spl+\Qv_\spl)$ is the coset of $\ov_\ell=\half(\ev_1+\dots+\ev_\ell)$ and the respective element of $T_\spl$ representing the non-unity component of $\PSO_{\ell,\ell}(\RR)$ is
\begin{equation*}
t=\exp\pi\ii\ov_\ell=\pm\diag(\underbrace{\ii,\dots,\ii}_\ell,\underbrace{-\ii,\dots,-\ii}_\ell).
\end{equation*}

If $\ell$ is even, then $2\Xvt_\spl\subseteq\Qv_\spl$ and $\pi_0\PSO_{\ell,\ell}(\RR)$ is of order~4. The two generators of order~2 of the 4-group $\Xv_\spl/(2\Xvt_\spl+\Qv_\spl)$ are the cosets of $\ov_\ell$ and of~$\ev_1$. The three non-unity components of $\PSO_{\ell,\ell}(\RR)$ are represented by the following elements of~$T_\spl$:
\begin{align*}
t_1&=\exp\pi\ii\ev_1=\pm\diag(-1,1,\dots,1,-1),\\
t_2&=\exp\pi\ii\ov_\ell=\pm\diag(\ii,\dots,\ii,-\ii,\dots,-\ii),\\
t_3&=t_1t_2=\pm\diag(-\ii,\ii,\dots,\ii,-\ii,\dots,-\ii,\ii).
\end{align*}

\subsection{}

Now let $G$ be an exceptional simple algebraic group. If $G$ is of type $\G_2$, $\F_4$, or $\E_8$, then $G(\CC)$ is always simply connected and therefore $G(\RR)$ is connected by the Cartan theorem. If $G$ is of type $\E_6$, then it is either simply connected or adjoint, and in the latter case $\pi_1G$ is of order~3. In both cases $G(\RR)$ is connected. The only remaining case is $\E_7$.

Let $G$ be a simple algebraic group of type $\E_7$ defined over $\RR$. Then $G$ is either simply connected or adjoint, and in the latter case $\pi_1G$ is of order~2. In the simply connected case, $G(\RR)$ is connected by the Cartan theorem. Consider the adjoint case.

If $G$ is anisotropic, then $G(\RR)$ is connected. It remains to consider the non-compact real forms of~$\E_7$. We use the notation of \cite[Tables]{Lie-3}.

The group $G$ contains a subgroup locally isomorphic to $\SL_8$ and sharing the maximal torus $T$ with~$G$. Let $\e_1,\dots,\e_8$ denote the weights of the tautological 8-dimensional representation of~$\SL_8$, which are obtained from the weights for $\GL_8$ (see Subsection~\ref{GL}) by restriction. They are bound by a unique linear dependence $\e_1+\dots+\e_8=0$. The coweights $\ev_1,\dots,\ev_8$ obtained from the respective coweights for $\GL_8$ by projection are bound by a similar linear dependence.

The coroots and cocharacters of $T$ are expressed in terms of $\ev_1,\dots,\ev_8$. The coroots of $G$ are $\ev_i-\ev_j$ and $\ev_i+\ev_j+\ev_k+\ev_l$, where the indices $i,j,k,l$ are pairwise distinct. The coroot lattice is
\begin{equation*}
\Qv=\{k_1\ev_1+\dots+k_8\ev_8\mid k_i\in\half\ZZ,\ k_i-k_j\in\ZZ,\ k_1+\dots+k_8=0\}.
\end{equation*}
A base of $\Qv$ is composed by the simple coroots
\begin{equation*}
\av_i=\ev_i-\ev_{i+1} \quad (i=1,\dots,6) \quad\text{and}\quad
\av_7=\ev_5+\ev_6+\ev_7+\ev_8.
\end{equation*}
The cocharacter lattice $\Xv$ is the coweight lattice of~$G$. It is spanned by coweights $\ev_i+\ev_j$. A base of $\Xv$ is composed by the fundamental coweights
\begin{equation*}
\ov_i=\ev_1+\dots+\ev_i+\min\{i,8-i\}\cdot\ev_8 \quad (i=1,\dots,6) \quad\text{and}\quad
\ov_7=2\ev_8.
\end{equation*}

If $G=\EV$ is the split real form of~$\E_7$, then $\Xv_\spl=\Xvt_\spl=\Xv$ and $\Qv_\spl=\Qv$. Hence $\pi_0G(\RR)\simeq\Xv/\Qv$ is of order~2. The only nonzero element of $\Xv/\Qv$ is the coset of $\ov_1$ and the respective element of $T_\spl=T$ representing the non-unity connected component $G(\RR)\setminus G(\RR)^\circ$ is $t=\exp\pi\ii\ov_1$.

If $G=\EVI$ is the quaternionic real form of~$\E_7$, then $\Xv_\spl$ is spanned by $\ov_2,\ov_4,\ov_5,\ov_6$; see \cite[Table~4]{Lie-3}. All these fundamental coweights belong to~$\Qv$. Hence $\Xv_\spl=\Qv_\spl$ and therefore $G(\RR)$ is connected.

If $G=\EVII$ is the Hermitian real form of~$\E_7$, then $\Xv_\spl$ is spanned by $\ov_1,\ov_2,\ov_6$ and $\ttt_\cmp$ is spanned by $\av_3,\av_4,\av_5,\av_7$; see \cite[Table~4]{Lie-3}. The lattice $\Qv_\spl$ is then spanned by $2\ov_1,\ov_2,\ov_6$ and $\Xvt_\spl$ is spanned by $\ov_1$, $\ov_2$, $\ov_6$, $(\ov_3)_\spl$, $(\ov_4)_\spl$, $(\ov_5)_\spl$, $(\ov_7)_\spl$, where
\begin{align*}
(\ov_3)_\spl &= \ov_3-\av_3-\av_4-\tfrac12\av_5-\tfrac12\av_7 = \ev_1+\ev_2-\tfrac12\ev_7+\tfrac52\ev_8 = \ov_2+\tfrac12\ov_6, \\
(\ov_4)_\spl &= \ov_4-\av_3-2\av_4-\av_5-\av_7 = \ev_1+\ev_2-\ev_7+3\ev_8 = \ov_2+\ov_6, \\
(\ov_5)_\spl &= \ov_5-\tfrac12\av_3-\av_4-\av_5-\tfrac12\av_7 = \tfrac12\ev_1+\tfrac12\ev_2-\ev_7+2\ev_8 = \tfrac12\ov_2+\ov_6, \\
(\ov_7)_\spl &= \ov_7-\tfrac12\av_3-\av_4-\tfrac12\av_5-\av_7 = \tfrac12(\ev_1+\ev_2-\ev_7+3\ev_8) = \tfrac12(\ov_2+\ov_6).
\end{align*}
It follows that $2\Xvt_\spl\subset\Qv_\spl$ and $\pi_0G(\RR)\simeq\Xv_\spl/\Qv_\spl$ is of order~2. The only nonzero element of $\Xv_\spl/\Qv_\spl$ is the coset of $\ov_1$ and the respective element representing the non-unity connected component of $G(\RR)$ is $t=\exp\pi\ii\ov_1$.

\begin{remark}
The cardinalities of $\pi_0G(\RR)$ for real forms of the adjoint group $G$ of type $\E_7$ were tabulated in \cite[Table~5]{Adams-Taibi}.
\end{remark}

\end{document}